\documentclass[12pt,a4paper,oneside,onecolumn]{article}
\usepackage[]{latexsym}

\newcommand{\nit}{{\rm I\!N}}

\newcommand{\rit}{{\rm I\!R}}
\newcommand{\cqfd}{\hspace*{\fill}\rule{2mm}{2mm}} 
\newcommand{\pre}{{\bf Proof.\ }}

\newcommand{\vu}{{\bf u}}
\newcommand{\vv}{{\bf v}}

\newcommand{\vA}{{\bf A}}
\newcommand{\vB}{{\bf B}}
\newcommand{\vV}{{\bf V}}
\newcommand{\vW}{{\bf W}}

\newcommand{\cs}{{\mathcal S}}

\newcommand{\ca}{{\mathcal A}}
\newcommand{\cb}{{\mathcal B}}
\newcommand{\cc}{{\mathcal C}}
\newcommand{\cd}{{\mathcal D}}
\newcommand{\ct}{{\mathcal T}}

\newcommand{\cp}{{\mathcal P}}
\newcommand{\cl}{{\mathcal L}}

\newcommand{\tr}{{\rm tr}} 
\newcommand{\pr}{{\rm pr}} 
\newcommand{\1}{{\bf 1}} 
\newcommand{\V}{{\rm Vect}}

\def\bvec#1{\overrightarrow{#1}} 

\newtheorem{theo}{Theorem}[section]
\newtheorem{prop}[theo]{Proposition}
\newtheorem{lema}[theo]{Lemma}
\newtheorem{coro}[theo]{Corollary}
\newtheorem{rema}[theo]{Remark}

\title{Affine representations of Lie algebras and geometric interpretation in the case of smooth manifolds}
\author{Sarah Hansoul \& Pierre Lecomte\footnote{The first author thanks the Fonds National de la Recherche Scientifique de Belgique for her Fellowship}\\
{\small Institut de Math\'ematique, Universit\'e de Li\`ege}\\
{\small Grande Traverse, 12 (B37), B-4000 Li\`ege, Belgium.}\\
{\small e-mail: s.hansoul@ulg.ac.be, plecomte@ulg.ac.be}}
\date{}
\begin{document}
\maketitle
\abstract{\small In order to unsderstand the structure of the cohomologies involved in the study of projectively equivariant quantization, we introduce a notion of affine representation of a Lie algebra.
We show how it is related to linear representations and $1$-cohomology classes of the algebra.
We classify the affine representations of the Lie algebra of vector fields associated to its action on symmetric tensor fields of type $(^1_2)$. Among them, we recover the space of symmetric linear connections and the space of projective structures. We compute some of the associated cohomologies.}

\section*{}

{\bf Key words.} Cohomology of Lie algebra, affine space, representation, differential manifold.

\noindent
{\bf Mathematics Subject Classifications (2000)}. 17B56, 17B66, 53C99.

\section{Introduction}

Cohomology classes of the Chevalley-Eilenberg cohomology of a representation of a Lie algebra are affine spaces in a natural way, modeled on the space of coboundaries. They bear an additionnal structure that comes from the action of the Lie algebra on the space of cochains. For $1$-cohomology classes, that structure is a special case of what we call an {\it affine representation} of the Lie algebra.

Differential geometry leads to natural examples of such representations. For instance, the set of linear connections of a smooth manifold $M$ is an affine representation of its Lie algebra of vector fields. This is also the case of the set $\cp_M$ of projective structures of the manifold, that is a quotient of its space $\cc_M$ of torsion free linear connections. Both of these examples play a role in the study of the so-called projectively equivariant quantizations introduced in \cite{10}. There, it is explained how the existence and uniqueness problem for these quantizations may be formulated in terms of cohomologies defined on some spaces of polynomials on $\cc_M$ and $\cp_M$. 

In that paper, however, the definitions of these cohomologies were only alluded to. Moreover, the notion of polynomials as well as that of the action of the Lie algebra of vector fields on these are not so obvious as the affine spaces involved are infinite dimensional. One of our purposes here is to provide a neat algebraic setting in which these notions can be easily handled.
 
Martin Bordemann has recently shown the existence of projectively equivariant quantization \cite{1}. His method makes no use of the cohomological approach. This approach is still  relevant to that question however. It remains indeed to classify the various equivariant quantizations, what seems to be difficult to perform using bordemann's framework. Moreover, the true origin of the so-called resonnant values \cite{6} still has to be found for arbitrary manifolds (they are obtained ``by force'' in Bordemann's paper while, on $\rit^m$, their cohomological origin is obvious \cite{11}).

After having defined the notion of affine representation of a Lie algebra $L$, we show that it induces a representation on a vector space $\vV$ of $L$ and a $1$-cohomology class $c$ of that representation. It is then proven that both the representation and the class determine the affine representation up to equivalence. It appears that $c$ is itself an affine representation of $L$ in a natural way. If $H^0(L,\vV)=0$, it is even equivalent to the affine representation that has induced it so that in this case, $H^1(L,\vV)$ determines completely the set of equivalence classes of affine representations of $L$ inducing the representation $\vV$. 

For instance, the representation of the Lie algebra $\V(M)$ of vector fields induced by the affine representation $\cc_M$ is the space $\cs^1_2(M)$ of symmetric tensor fields of $M$ of type $(^1_2)$ equipped with the Lie derivative while that corresponding to $\cp_M$ is the subspace of $\cs^1_2(M)$ consisting of these elements whose contraction is vanishing.

We then recall the definition of the space $\cp(\ca,\vB)$ of polynomials on an affine space $\ca$ valued in some vector space $\vB$, taking care of the fact that $\ca$ or $\vB$ could possibly be infinite dimensional. We show that, when $\ca$ is an affine representation of $L$ and $\vB$ is a linear representation, then $\cp(\ca,\vB)$ is a representation of $L$ in a natural way. In this case, the short exact sequences
\[
0\to\cp^{k-1}(\ca,\vB)\to\cp^k(\ca,\vB)\to S^k(\bvec{\ca},\vB)\to 0
\]
induced by the filtration of $\cp(\ca,\vB)$ by the degree $k$ of polynomials are, in general, not split and closely related to the cohomology class induced by $\ca$.

As an application of the above considerations, we classify the affine representations of $\V(M)$ inducing $\cs^1_2(M)$ and we compute the cohomology spaces
\[
H^i(\V(M),\cp^k_{diff}(\cc_M,\cs^1_2(M))), i=0,1,
\]
where the subscript $diff$ means that we restrict ourselves to polynomials whose components relative to any base point $\nabla\in\cc_M$ are multidifferential operators.
This computation is quite long and hard and we only sketch its main key steps. It uses the above short exact sequences.
\tableofcontents

\section{Algebraic preliminaries}

\subsection{Affine representations of a Lie algebra}\label{secaffrep}
Let $\ca$ be an affine space modeled on a vector space $\vA$.
An {\it affine repesentation} of a Lie algebra $L$ on $\ca$ is a map
\[
(x,a)\in L\times\ca\mapsto x.a\in\vA
\]
such that each map $a\mapsto x.a$, $x\in L$, is affine from $\ca$ into $\vA$ (equipped with its canonical affine structure) and such that, for all $x,y\in L$ and $a\in\ca$,
\begin{equation}\label{affrep}
\bvec{x}.(y.a)-\bvec{y}.(x.a)-[x,y].a=0
\end{equation}
(where $\overrightarrow{x}$ denotes the linear map associated to $a\mapsto x.a$)(\footnote{For simplicity, we will also say that $\ca$, equipped with such a map, is an affine representation of $L$.}).

Given affine representations of $L$ on affine spaces $\ca$ and $\cb$, an {\it intertwinning map}, also said {\it affine equivariant map}, from $\ca$ into $\cb$ is an affine map $f:\ca\to \cb$ such that
\[
x.f(a)=\bvec{f}(x.a)
\]
for all $x\in L$ and $a\in \ca$. We will say that the representations $\ca$ and $\cb$ are {\it isomorphic} or {\it equivalent} if there is a bijective intertwinning map from $\ca$ into $\cb$.
We will denote by $\ca ff_L(\ca,\cb)$ the set of intertwinning maps from $\ca$ into $\cb$.

Let us replace $a$ by $a_0+\vu$ in (\ref{affrep}), where $a_0\in\ca$ is fixed and $\vu\in\vA$ is arbitrary. Substracting member by member (\ref{affrep}) from the equality so obtained, we first get the relation
\[
\bvec{x}.(\bvec{y}.\vu)-\bvec{y}.(\bvec{x}.\vu)-\bvec{[x,y]}.\vu=0.
\]
It shows that $x\in\vA\mapsto\bvec{x}\in gl(\vA)$ is a (linear) representation $\rho_\ca$ of $L$.
We then obtain
\[
\rho_\ca(x)(y.a_0)-\rho_\ca(y)(x.a_0)-[x,y].a_0=0
\]
which means that $\gamma_{a_0}:x\mapsto x.a_0$ is a $1$-cocycle of $L$ valued in $\vA$. In addition, since
\[
\gamma_a=\gamma_{a_0}+\partial(a-a_0),
\]
where $\partial$ denotes the Chevalley-Eilenberg coboundary operator, we see that the cohomology class $c_\ca$ of $\gamma_{a_0}$ does not depend on $a_0$.

Hence, an affine representation $\ca$ produces a representation $\rho_\ca$ of $L$ on $\vA$ and a class $c_\ca\in H^1(L,\vA)$ of the cohomology of that representation. We call them the representation and the cohomology class {\it induced} by the affine representation $\ca$.

Conversely, let $\rho$ be a representation of $L$ on $\vA$ and let $c$ be an element of its first cohomology space. Choose $\gamma_0\in c$. Viewing $\vA$ as an affine space in the canonical way, it is easy to see that
\[
(x,\vu)\in L\times\vA\mapsto\gamma_0+\rho(x)\vu\in\vA
\]
is an affine representation of $L$ on $\vA$ inducing $\rho$ and $c$. Moreover, if $\rho=\rho_\ca$, $c=c_\ca$ and $\gamma_0=\gamma_{a_0}$, then this affine representation is equivalent to $\ca$, a bijective intertwinning map from $\ca$ into $\vA$ being
\[
\varphi:a\in\ca\mapsto a-a_0\in\vA.
\]
Indeed, $\bvec{\varphi}$ is the identity map of $\vA$ and
\[
\bvec{\varphi}(x.a)=x.a=x.(a_0+(a-a_0))=\gamma_0(x)+\rho_\ca(x)(a-a_0)=x.\varphi(a).
\]

On the other hand, the following proposition is easily checked.

\begin{prop}\label{intertwinning}
Let $\ca$ and $\cb$ be affine representations of $L$. An affine map $f:\ca\to\cb$ is an intertwinning map from $\ca$ into $\cb$ if and only if $\bvec{f}$ intertwinnes $\rho_\ca$ and $\rho_\cb$ and $c_\cb=\bvec{f}_\sharp c_\ca$.
\end{prop}

Let us say that $(\rho,c)$ and $(\rho',c')$, where $\rho, \rho'$ are representations of $L$ and $c,c'$ are $1$-cocycles of their respective cohomologies, are {\it equivalent} if there is a bijective intertwinning map between $\rho$ and $\rho'$ that transforms $c$ into $c'$. Then, the above observations could be summarized as follows.

\begin{prop}
The map $\ca\mapsto (\rho_\ca,c_\ca)$ induces a bijective correspondance between the set of equivalence classes of affine representations of $L$ and the set of equivalence classes of couples $(\rho,c)$ consisting of a representation $\rho$ of $L$ and a class $c$ of degree $1$ of its cohomology.
\end{prop}

\begin{coro}
Let $\ca$ and $\cb$ be affine representations of $L$. If the number $\lambda$ is not equal to $zero$ and if $\rho_\cb=\lambda \rho_\ca, c_\cb=\lambda c_\ca$, then $\ca$ and $\cb$ are equivalent.
\end{coro}

Recall that the {\it direct sum} $\ca_1\oplus\ca_2$ of the affine spaces $\ca_1, \ca_2$ modeled respectively on vector spaces $\vA_1,\vA_2$ is the space $\ca_1\times\ca_2$ on which the space $\vA_1\oplus\vA_2$ acts by the translations
\[
(a_1,a_2)+(u_1,u_2)=(a_1+u_1,a_2+u_2).
\]
If $\ca_1, \ca_2$ are affine representations of $L$, then $\ca_1\oplus\ca_2$ is also an affine representation of $L$, for the action
\[
x.(a_1,a_2)=(x.a_1,x.a_2).
\]

\begin{coro}
Let $\ca$ be an affine repesentation of $L$. Assume that $\vA=\vA_1\oplus\vA_2$ where the $\vA_i$ are stable subspaces for $\rho_\ca$ and let $\pi_i:\vA\to\vA_i, i=1,2,$ denote the corresponding projectors. For $i=1,2$, let also $\ca_i$ be an affine representation inducing the representation $\rho_{\ca_i}=\pi_i\circ\rho_\ca$ and its cohomology class $c_{\ca_i}=\pi_{i\sharp}c_\ca$. Then $\ca$ is equivalent to the direct sum of affine representations $\ca_1\oplus\ca_2$.
\end{coro}

The two above corollaries are straightforward. They can be used to classify the affine representations of $L$ inducing a given representation.

\medskip
Let us have a closer look at the cohomology class $c_\ca$ induced by an affine representation $\ca$ of $L$. It is an affine space modeled on the vector space $\partial\vA$, the space of $1$-coboundaries of $\rho_\ca$. It carries itself an affine representation of $L$, given by
\[
(x,\gamma)\in L\times c_\ca\mapsto\partial(\gamma(x))\in\partial\vA.
\]

 \begin{prop}
The map $\Phi:a\mapsto \gamma_a$ is an intertwinning map from the affine representation $\ca$ into the affine representation $c_\ca$ of $L$. It is onto and $\Phi(a')=\Phi(a)$ if and only if $\partial(a'-a)=0$. In particular, $c_\ca$ is isomorphic to $\ca/H^0(L,\vA)$.
\end{prop}
\pre This is straightforward.$\cqfd$

\medskip
Let $\rho$ be a representation of $L$ on the vector space $\vV$. Suppose that $H^0(L,\vV)=0$. This means that for $\vv\in\vV$, $\partial\vv=0$ implies $\vv=0$. Any cohomology class $c\in H^1(L,\vV)$ can then be viewed as an affine space modeled on $\vV$, $\vv$ acting on $c$ by the translation
\[
\gamma\mapsto \gamma+\partial\vv.
\]
With that identification, the canonical affine representation of $L$ defined above on $c$ is given by the evaluation:
\[
x.\gamma=\gamma(x)
\]
The above Proposition (\ref{intertwinning}) can then be specialized as follows.
\begin{prop}\label{specialcase}
Let $\rho$ be a representation of $L$ on the vector space $\vV$ such that $H^0(L,\vV)=0$. Let $c,c'\in H^1(L,\vV)$. A map $f:c\to c'$ is an affine equivariant map if and only if it is of the form
\[
\gamma\in c\mapsto T\circ\gamma\in c'
\]
for some $T\in H^0(L,Hom(\vV,\vV))$ such that $c'=T_\sharp c$.
\end{prop}

Hence, the elements of $H^1(L,\vV)$ are affine representations of $L$ inducing $\rho$ and each affine representation of $L$ inducing $\rho$ is equivalent to some element of $H^1(L,\vV)$. 
Combining the above propositions and corollaries, we then get the following proposition.
\begin{prop}
Let $\rho$ be a representation of $L$ on the vector space $\vV$. Suppose that $H^0(L,\vV)=0$ and that $\dim H^0(L,Hom (\vV,\vV))=1$. If $\dim H^1(L,\vV)=s$, then there exists exactly $2^s$ equivalence classes of affine representations of $L$ inducing the representation $\rho$.
\end{prop}

\subsection{Polynomials on an affine space}\label{afrep}
We recall here some useful facts about polynomials defined on affine spaces, possibly infinite dimensional. The reader is referred to \cite{2} for further details about that notion.

Let $\vA$ be a vector space. For any map $f$ defined on $\vA^i$, $\hat{f}$ denotes the composition $f\circ\Delta$ of $f$ and the diagonal map $\Delta:\vu\in\vA\mapsto (\vu,\ldots,\vu)\in \vA^i$.
If $f$ takes its values in a vector space $\vB$ and if it is $i$-linear and symmetric, then it is completely determined by $\hat{f}$ as stated in the following lemma.

\begin{lema}
If $f,g : \vA^i\to\vB$ are $i$-linear symmetric maps then $f=g$ if and only if $\hat{f}=\hat{g}$.
\end{lema}
\pre See \cite{2}.$\cqfd$

\medskip
Let $\ca$ be an affine space modeled on the vector space $\vA$ and $\vB$ be a vector space.

A mapping $p: \ca\to\vB$ is a {\it polynomial of degree $\leq k$} if, for some $a_0\in\ca$, 
\begin{equation}\label{expressionp}
p(a_0+\vu)=\sum_{i=0}^k\frac 1 {i!}\hat{p}_i(\vu),\ \forall \vu\in\vA,
\end{equation}
for some symmetric $i$-linear maps $p_i$ from $\vA^i$ into $\vB$. It follows from the above Lemma that these are uniquely determined by $p$ and the base point $a_0$ chosen in $\ca$. We call them the {\it components of $p$ relative to the base point $a_0$}.

The above definition is independant of that point. Indeed

\begin{lema}\label{changeexpressionp}
If $a'_0\in\ca$, assuming that $p$ is defined by {\rm (\ref{expressionp})},
\[
p(a'_0+\vu)=\sum_{i=0}^k\frac 1 {i!}\hat{q}_i(\vu),\ \forall \vu\in\vA,
\]
where
\[
q_j(\vu_1,\ldots,\vu_j)=\sum_{i\geq j}\frac 1 {(i-j)!}p_i(\underbrace{a'_0-a_0,\ldots,a'_0-a_0}_{i-j},\vu_1,\ldots,\vu_j)
\]
for all $j\in\{0,\ldots,k\}$.
\end{lema}
\pre This is straightforward, using the previous lemma.$\cqfd$

\medskip
We denote by $\cp^k(\ca,\vB)$ the set of polynomials from $\ca$ into $\vB$. It is a vector space. The sequence of spaces $\cp^k(\ca,\vB)$ is an increasing filtration. We denote by $\cp(\ca,\vB)$ its limit. It is the space of {\it $\vB$-valued polynomials defined on $\ca$}.

It follows from the previous Lemma that the vector space
\[
\cp^k(\ca,\vB)/\cp^{k-1}(\ca,\vB)
\]
is canonically isomorphic to the space $S^k(\vA,\vB)$ of $k$-linear maps from $\vA^k$ into $\vB$. We thus have a short exact sequence of vector spaces
\begin{equation}\label{shortsequence}
  0\to \cp^{k-1}(\ca,\vB)\rightarrow \cp^k(\ca,\vB)\rightarrow S^k(\vA,\vB)\to 0.
\end{equation}

By definition, any choice of a base point $a_0$ in $\ca$ defines a splitting of that sequence. However, it has no canonical splitting. We will see soon circumstances turning it into a sequence of Lie algebra representations. It then will no longer automatically be split.

\subsection{Cohomologies related to  affine representations}
Let $\ca$ be an affine representation of $L$, modeled over the vector space $\vA$, and let $\vW$ be a linear representation of $L$. We will see that the space $\cp(\ca,\vW)$ is endowed with a natural representation of $L$. We define it using a base point $a_0\in \ca$ but it does not depend on any particular choice of that point. Recall that the space $S^k(\vA,\vW)$ is equipped with the tensor product representation of copies of $\vA$ and $\vW$. It is given by, for $x\in L$ and any $\vu_i\in \vA$,
\[
(x.t)(\vu_1,\ldots,\vu_k)=x.(t(\vu_1,\ldots,\vu_k))-\sum_it(\vu_1,\ldots,x.\vu_i,\ldots,\vu_k).
\]
Let $x\in L$ and let $p\in\cp^k(\ca,\vW)$ have components $p_i$ relative to some $a_0\in \ca$. Then, by definition, the components relative to $a_0$ of the polynomial $x.p\in \cp^k(\ca,\vW)$ are given by
\begin{equation}\label{action}
(x.p)_i(\vu_1,\ldots,\vu_i)=(x.p_i)(\vu_1,\ldots,\vu_i)-p_{i+1}(x.a_0,\vu_1,\ldots,\vu_i).
\end{equation}

\begin{prop}
a) The polynomial $x.p$ does not depend on the base point $a_0$ of the affine space $\ca$.\\
b) The map $(x,p)\in L\times\cp^k(\ca,\vW)\mapsto x.p\in\cp^k(\ca,\vW)$ defines a representation of $L$ on $\cp^k(\ca,\vW)$.\\
\end{prop}
\pre This is easy to check by direct computation, using Lemma \ref{changeexpressionp}.$\cqfd$

\paragraph{Associated short exact sequences}
When the spaces $\cp^l(\ca,\vW)$ are equipped with the above action of $L$, then the sequence (\ref{shortsequence}) becomes a short exact sequence of representations of $L$. We will call it the {\it $k$-th exact sequence associated to the affine space $\ca$ and to $\vW$}.

\begin{prop}
 Let $\ca$ and $\ca'$ be affine representations of $L$ modeled on the same vector space $\vA$. Assume that $\rho_{\ca'}=\rho_\ca$. The $1$-th exact sequences associated to $\ca$ and $\ca'$ and to $\vA$ are isomorphic if and only if $c_{\ca'}=c_\ca.$ If $\vA$ is finite dimensional, then the same is true about the $1$-th sequences associated to the trivial representation $\rit$.
\end{prop}
\pre Recall that the isomorphism class of the sequence (\ref{shortsequence}) is an element $\alpha^k_{\ca,\vW}$ of $H^1(L,Hom(S^k(\vA,\vW),\cp^{k-1}(\ca,\vW)))$. It is constructed as follows. Let the section $\tau:S^k(\vA,\vW)\to\cp^k(\ca,\vW)$ of the projection $\cp^k(\ca,\vW)\to S^k(\vA,\vW)$ be defined by
\[
\tau(t)(a_0+\partial \vu)=\frac 1 {k!}\hat{t}(\vu).
\]
The coboundary $\partial\tau$ takes its values in $Hom(S^k(\vA,\vW),\cp^{k-1}(\ca,\vW))$. For $x\in L$, $(\partial\tau)_x$ maps indeed $t\in S^k(\vA,\vW)$ onto the polynomial of degree $k-1$
\[
a_0+\partial \vu\mapsto -\frac1{(k-1)!}t(x.a_0,\vu,\ldots,\vu).
\]
The class $\alpha^k_{\ca,\vW}$ is then the cohomology class of the cocycle $x\mapsto (\partial\tau)_x$.
(See \cite{4} for further details about the classification of short exact sequences.)

It follows from that construction that the classes $\alpha^1_{\ca,\vA}$ and $\alpha^1_{\ca',\vA}$ are equal if and only if there is a linear map $\theta:gl(\vA)\to \vA$ such that
\[
t(x.a'_0)=t(x.a_0)+x.\theta(t)-\theta(x.t)
\]
for all $x\in L$ and all $t\in gl(\vA)$, where $a_0\in \ca$ and $a'_0\in \ca'$.
Thus, if they are equal, then, $id_\vA$ denoting the identity from $\vA$ into itself, we get
\[
x.a'_0=x.a_0+x.\theta(id_\vA), \ \forall x\in L,
\]
so that $c_{\ca'}=c_\ca$. Conversely, if $c_{\ca'}=c_\ca$, then $x.a'_0=x.a_0+x.\vu$ for some $\vu\in \vA$ and we can take $\theta(t)=t(\vu)$ to show that  $\alpha^1_{\ca',\vA}=\alpha^1_{\ca,\vA}$.

In a similar fashion, we see that the classes $\alpha^1_{\ca,\rit}$ and $\alpha^1_{\ca',\rit}$ are equal if and only if there is a linear map $\mu$ from the dual $\vA^*$ of $\vA$ into $\rit$ such that
\[
\xi(x.a'_0)=\xi(x.a_0)-\mu(x.\xi),
\]
for all $x\in L$ and all $\xi\in \vA^*$. Suppose that this latter condition holds true and that $\vA$ is finite dimensional. Then $\vA^{**}\cong \vA$ in a natural way so that, viewing $\mu$ as an element of $\vA$, the above relation reads $x.a'_0=x.a_0+x.\mu$. Hence $\alpha^1_{\ca',\rit}=\alpha^1_{\ca,\rit}$ implies $c_{\ca'}=c_\ca$. The converse is clear. $\cqfd$

\medskip
There is a corresponding long exact sequence of cohomology spaces associated to the sequence (\ref{shortsequence}). It can be pictured as
\begin{equation}\label{longsequence}
\begin{array}{ccc}
H(L,\cp^{k-1}(\ca,\vW))&&\\
&\searrow&\\
\chi\uparrow&&H(L,\cp^k(\ca,\vW))\\
&\swarrow&\\
H(L,S^k(\vA,\vW))&&
\end{array}
\end{equation}
where $\chi$ is the connecting homomorphism (see \cite{7}). The degree of the latter is $1$. As easily seen, with the notations of the above proof, it maps the cohomology class of any $p-$cocycle
\[
t:(x_1,\ldots,x_p)\in L^p\mapsto t_{x_1,\ldots,x_p}\in S^k(\vA,\vW)
\] 
onto that of the $(p+1)$-cocycle 
\[
(x_0,\ldots,x_p)\in L^{p+1}\mapsto \tau(t^\chi_{x_0,\ldots,x_p})\in\cp^{k-1}(\ca,\vW)
\]
where
\[
t^\chi_{x_0,\ldots,x_p}:(\vu_1,\ldots,\vu_{k-1})\mapsto \sum_{i=0}^p(-1)^{i+1}t_{x_0,\ldots\hat{i}\ldots,x_p}(x_i.a_0,\vu_1,\ldots,\vu_{k-1})
\]
($\hat{i}$ denotes the omission of $x_i$).

\begin{prop}
Let $\ca,\cb$ be affine representations of $L$.  For each $f\in \ca ff_L(\ca,\cb)$,
\[
f^*:\cp(\cb,\vW)\to\cp(\ca,\vW):q\mapsto q\circ f
\]
is an intertwinning map.
\end{prop}
\pre As easily seen, the components of $p=f^*q$ relative to some base point $a_0\in \ca$ are given by
\[
p_i=\bvec{f}^*q_i:(\vu_1,\ldots,\vu_i)\mapsto q_i(\bvec{f}(\vu_1),\ldots,\bvec{f}(\vu_i))
\]
where the $q_i$'s are the components of $q$ relative to $f(a_0)\in \cb$. Using Proposition \ref{intertwinning} and the definition (\ref{action}) of the action of $L$ on polynomials, it is then a matter of elementary computations to conclude. $\cqfd$

\section{Application to smooth manifolds}
Let $M$ be a smooth manifold, connected, Hausdorff and second countable.
We set $\dim M=m$, and we assume that $m>1$.

The space $\cs^1_2(M)$ of symmetric tensor fields of $M$ of type $(^1_2)$, equipped with the Lie derivative, is a representation of the Lie algebra $\V(M)$ of vector fields of $M$.

We first classify the affine representations of $\V(M)$ that induce that representation.
As $H^0(\V(M),\cs^1_2(M))=0$, it follows from the end of subsection \ref{secaffrep} that each of these is equivalent to the affine representation defined on some element of $H^1(\V(M),\cs^1_2(M))$. We will see that the dimension of that space is $2$ and that there are $4$ equivalence classes of affine representations inducing $\cs^1_2(M)$.

More precisely, we will see that $\cs^1_2(M)$ splits as the direct sum of two $\V(M)$-modules. One of them is isomorphic to the space $\Omega_1(M)$ of smooth $1$-forms of $M$ equipped with the Lie derivative and its first cohomology space is one dimensional. The first cohomology space of the second is also one dimensional and admits a generator that, as an affine space, is canonically isomorphic to the space $\cp_M$ of projective structures of $M$. The space  $H^1(\V(M),\cs^1_2(M))$ is generated by the classes induced by these of the submodules. It has an element which is isomorphic to the affine space $\cc_M$ of linear torsion free connections of $M$.

The section ends with the computations of the first cohomology spaces of the subspace of $\cp^k(\cc_M,\cs^1_2(M))$ whose elements are multidifferential operators.
These are quite difficult and long. We will restrict ourselves to the description of the strategy of computation, giving only details for the key arguments.

\subsection{The representation $\cs^1_2(M)$}

Let us introduce the following linear maps
\[
\tr: \sum_{ijk}S^k_{ij}\epsilon^i\otimes\epsilon^j\otimes e_k\in \vee^2\rit^{m*}\otimes \rit^m \mapsto \sum_{ij}S^j_{ij}\epsilon^i\in\rit^{m*},
\]
 which we will call the \textit{trace}, and
\[
\alpha\in\rit^{m*}\mapsto \alpha\1\in\vee^2\rit^{m*}\otimes \rit^m,
\]
where
\[
\alpha\1:(\vu,\vv)\mapsto \alpha(\vu)\vv+\alpha(\vv)\vu.
\]
(We denote by $\1$ the idendity matrix.)
Both maps are $Gl(m,\rit)$-equivariant.

Observe that $\tr(\alpha\1)=(m+1)\alpha$ for all $\alpha$. In particular, for any $S$,
\[
\pr(S)=S-\frac{1}{m+1}\tr(S)\1
\]
has a vanishing trace. It will be convenient to denote also $\pr(S)$ by $\overline{S}$.
\begin{lema}
The representation $\vee^2\rit^{m*}\otimes\rit^m$ of $Gl(m,\rit)$ splits as the direct sum of two irreducible representations namely the kernel and the image of the linear equivariant map
\[
S\mapsto \tr(S)\1.
\]
\end{lema}
\pre This is straightforward using for instance the classical theory of invariant \cite{5}. $\cqfd$

\medskip
Observe that the image of the above map is just a copy of the irreducible representation $\rit^{m*}$ of $Gl(m,\rit)$.

\medskip
The space $\cs^1_2(M)$ is the space of sections of the bundle $\vee^2T^*M\otimes TM$. The latter is associated to the bundle of linear frames of $M$ and to the representation $\vee^2\rit^{m*}\otimes\rit^m$ of its structure group $Gl(m,\rit)$. It thus splits as a direct sum of two subbundles according to the decomposition of that representation into irreducible components. The maps induced on $\vee^2T^*M\otimes TM$ in a natural way by $\pr$ and $\tr$ will be denoted by the same notations. The restriction of $\tr$ to one of the subbundles is an isomorphism onto $T^*M$. Taking sections, we get a decomposition 
\begin{equation}\label{decomposition}
\cs^1_2(M)=\overline{\cs^1_2(M)}\oplus \tr(\cs^1_2(M))\1
\end{equation}
of $\cs^1_2(M)$ as a direct sum of representations of $\V(M)$ (as it induces no real confusion, we use also the notation $\1$ to denote the identity endomorphism of $TM$). Moreover, the representation $\tr(\cs^1_2(M))\1$ is isomorphic to that defined by the Lie derivative acting on the space $\Omega_1(M)$ of smooth $1$-forms of $M$.

As above, we will  keep the same notations to denote the maps induced by $\pr$ and $\tr$ on the space of sections of $\vee^2T^*M\otimes TM$ in the natural way.

\begin{prop}\label{invariant}
$H^0(\V(M),\cs^1_2(M))=0.$
\end{prop}
\pre It is indeed a matter of simple calculation to see, using local coordinates, that if $L_XS=0$ for all $X\in\V(M)$, then $S=0$.$\cqfd$

\begin{prop}\label{entrelacementS12}
A linear map
\[
T:\cs^1_2(M)\to\cs^1_2(M)
\]
is $\V(M)$-equivariant if and only if there exist real numbers $a$ and $b$ such that
\begin{equation}\label{entrelacement}
T(S)=a\overline{S}+b\, \tr(S)\1, \forall S\in \cs^1_2(M).
\end{equation}
\end{prop}
\pre It is clear that the maps of the form (\ref{entrelacement}) are equivariant.
Conversely, let $T: \cs^1_2(M)\to \cs^1_2(M)$ be equivariant. It follows immediately from Proposition 4 and Example 12 of \cite{3} that it is a differential operator (this uses the assumption $\dim M >1$).  In a local chart of $M$  of connected domain $U$, expressing the invariance with respect to constant and linear vector fields, one sees then that the local expression of $T$ is $gl(m,\rit)$-equivariant. Using the theory of invariants (\cite{5}, for instance), one then shows easily that $T_{|U}=a_U\pr_{|U}+b_U\tr_{|U}\1$ for some constants $a_U$ and $b_U$. It is clear that if connected domains of charts $U$ and $V$ have nonempty intersection, then $a_U=a_V$ and $b_U=b_V$. Hence the conclusion as $M$ is connected.  $\cqfd$

\medskip
It will be useful to denote by $\tr.\1$ the map $S\mapsto \tr(S)\1$.

\begin{coro}
The multiples of $\pr$ and $\tr$ are exactly the $\V(M)$-equivariant maps from $\cs^1_2(M)$ into $\overline{\cs^1_2(M)}$ and $\Omega_1(M)$ respectively.
\end{coro}
\pre This is straightforward. $\cqfd$

\subsection{$H^1(\V(M),\cs^1_2(M))$}

The space $\cc_M$ is an affine representation of $\V(M)$ in a natural way. The representation is defined by
\[
\cl: (X,\nabla)\in\V(M)\times\cc_M\mapsto L_X\nabla\in\cs^1_2(M)
\]
where
\[
L_X\nabla(Y,Z)=[X,\nabla_YZ]-\nabla_{[X,Y]}Z-\nabla_Y[X,Z],\forall Y,Z\in \V(M).
\]
Indeed, it is known \cite{9} that for each $\nabla\in\cc_M$,
\[
\cl_\nabla:X\mapsto L_X\nabla
\]
is a cocycle whose cohomology class is independant of $\nabla$.

The maps $\pr\circ\cl_\nabla$ and $\tr\circ\cl_\nabla$ are cocycles too. As $\overline{\cs^1_2}(M)\subset\cs^1_2(M)$, we can view $\pr\circ\cl_\nabla$ as valued in the latter, in which case we denote it $\cl^{pr}_\nabla$. We also set $\cl^{tr}_\nabla=(\tr\circ\cl_\nabla)\1$.

\begin{prop}\label{cohomologieS12}
For each $\nabla\in\cc_M$, the cohomology classes of $\cl^{pr}_\nabla$ and $\cl^{tr}_\nabla$ form a basis of the space $H^1(\V(M),\cs^1_2(M))$.
\end{prop}
\pre Since the maps $\pr$ and $\tr$ are $\V(M)$-equivariant, $\cl^{pr}_\nabla$ and $\cl^{tr}_\nabla$ are cocycles. They are not coboundaries because they are differential operators from $\V(M)$ into $\cs^1_2(M)$ of order $2$ while coboundaries are of order $1$. Their cohomology classes are linearly independant because they take their values into the supplementary submodules $\overline{\cs^1_2(M)}$ and $\tr(\cs^1_2(M))\1$ (see (\ref{decomposition})).

Each map $\gamma:\V(M)\to \cs^1_2(M)$ induces by contraction a linear map $^0\gamma:\V(M)\to\cd(\cs^2(M),\cs^1(M))$, $\cs^k(M)$ denoting the space of $k$-symmetric contravariant tensor fields of $M$ and $\cd$ denoting the space of differential operators. It is a cocycle or a coboundary if $\gamma$ is a cocycle or a coboundary respectively.

It has been shown in \cite{12} that
\[
H^1(\V(M),\cd(\cs^2(M),\cs^1(M)))
\]
is generated by the class $[^0\cl_\nabla]$ of $^0\cl_\nabla$. Moreover, it follows from \cite{10} that
\begin{equation}\label{prettr}
p\;^0\cl_\nabla+q\;^0\cl^{tr}_\nabla=\partial D_\nabla
\end{equation}
for some (nonzero) real numbers $p$ and $q$, where $D_\nabla\in\cd(\cs^2(M),\cs^1(M))$ maps $P$ onto the contraction of the covariant differential $\nabla P$ of $P$. Given a $1$-cocycle $\gamma$, one has 
\[
^0\gamma=a^0\cl_\nabla+\partial\ct
\]
for some number $a$ and some linear differential operator $\ct:\cs^2(M)\to\cs^1(M)$. It follows from the lemma below that $^0\gamma-a^0\cl_\nabla$ is of the form $b\partial D_\nabla+^0(\partial T)$ for some number $b$ and some $T\in\cs^1_2(M)$. Using (\ref{prettr}), the (obvious) relation
\[
\cl_\nabla=\cl_\nabla^{\pr}+\frac{1}{m+1}\cl_\nabla^{\tr}
\]
and the fact that the map $\gamma\mapsto ^0\gamma$ is one-to-one, it follows that $\gamma-\partial T$ is a linear combination of $\cl^{\pr}_\nabla$ and $\cl^{\tr}_\nabla$.$\cqfd$
\begin{lema}\label{odrezero}
Let $\ct\in\cd(\cs^p(M),\cs^q(M))$, where $p > q$, be given.
If the Lie derivatives $L_X\ct, X\in\V(M),$ are operators of order $0$, then $\ct$ is of the form
\[
P\mapsto k\,D^{p-q}_\nabla(P)+^0T(P)
\]
for some some $T\in\cs^1_2(M)$ and some $k\in\rit$.
\end{lema}
\pre Suppose that $\ct$ is of positive order. Its principal symbol $\sigma_\ct$ is thus $\V(M)$-invariant. Indeed $L_X\sigma_\ct=\sigma_{L_X\ct}=0$ because the Lie derivatives of $\ct$ are of order $0$. As easily seen $\sigma_\ct:T^*M\oplus\cs^p(M)\to\cs^q(M)$ is thus a constant multiple of
\[
(\eta,P)\mapsto \iota^{p-q}_\eta P.
\]
This is precisely the principal symbol of $D^{p-q}_\nabla$. $\cqfd$

\medskip
\begin{coro}
For each $\nabla\in\cc_M$, one has
\[
H^1(\V(M),\overline{\cs^1_2(M)})=\rit [\pr\circ\cl_\nabla]
\]
and
\[
H^1(\V(M),\Omega_1(M))=\rit [\tr\circ\cl_\nabla].
\]
\end{coro}
\pre Straightforward in view of (\ref{decomposition}). $\cqfd$

\begin{rema}
{\rm  The space $H^1(\V(M),\Omega_1(M))$ was already known. See \cite{15} for instance.}
\end{rema}

\begin{rema}
{\rm Proposition \ref{cohomologieS12} shows that $H^1(\V(M),\cs^1_2(M))$ is generated by the classes found in \cite{9}. We could also have computed this space using \cite{15} although it would have been quite difficult to extract from that reference the explicit forms of the cocycle generating the cohomology that we need in the sequel.}
\end{rema}

\subsection{Affine representations inducing $\cs^1_2(M)$}
 As they do not depend on the choice of $\nabla$, the cohomology classes of the cocycles $\cl_\nabla$, $\cl^{pr}_\nabla$ and $\cl^{tr}_\nabla$ will be denoted $\cl_M$, $\cl^{pr}_M$ and $\cl^{tr}_M$ respectively.
\begin{prop}
Each affine representation of $\V(M)$ inducing $\cs^1_2(M)$ is isomorphic to  either $0$, $\cl_M$, $\cl^{pr}_M$ or $\cl^{tr}_M$, and to only one of these affine representations.
\end{prop}
\pre This follows easily from Proposition \ref{specialcase} and Proposition \ref{entrelacementS12}. $\cqfd$

\begin{prop}
As affine representations of $\V(M)$, the affine spaces $\cc_M$ and $\cl_M$ are isomorphic: the map
\[
\nabla\mapsto \cl_\nabla
\]
is natural, bijective and belongs to $\ca ff_{\V(M)}(\cc_M,\cl_M)$.
\end{prop}
\pre The map defined above is obviously natural, affine and onto. The local expression 
\[ 
\cl_\nabla(X)^k_{ij}=\partial_{ij}X^k+\partial_i X^u\Gamma ^k_{uj}+\partial_j X^u\Gamma ^k_{iu}-\partial_u X^k\Gamma ^u_{ij}+X^u\partial_u\Gamma^k_{ij}
\]
shows that it is one to one. $\cqfd$

\medskip
The map $\pr$ commutes with the Lie derivatives. The affine representation of $\V(M)$ on $\cc_M$ can thus be pushed down on the quotient $\cc_M/\Omega_1(M)$ of $\cc_M$ by $\ker\pr\cong\Omega_1(M)$. This affine representations induces the $\V(M)$-submodule $\overline{\cs^1_2(M)}$ of $\cs^1_2(M)$.
By a theorem of Weyl \cite{8}, there is a canonical bijection between $\cc_M/\Omega_1(M)$ and the set $\cp_M$ of projective structures of $M$. The latter is thus an affine representation of $\V(M)$ in a natural way. The following two corollaries are then straightforward.
\begin{coro}
The map $\nabla\mapsto\pr\circ\cl_\nabla$ induces a natural intertwinning map between the space $\cp_M$ of projective structures on $M$ and the cohomology class $\pr_\sharp\cl_M\in H^1(\V(M),\overline{\cs^1_2(M)})$.
\end{coro}
\begin{coro}
The affine representations $\cp_M\times\Omega_1(M)$ and $\cl_M^\pr$ are equivalent.
\end{coro}

\begin{rema}{\rm In a similar way, one can also identify the class $\cl^{tr}_M$  with  $\overline{S^1_2(M)}\times\tr_\sharp\cl_M$. The space $\tr_\sharp\cl_M$ can be viewed as the quotient of $\cl_M$ by $\overline{\cs^1_2(M)}$. However, we have no geometric interpretation of that affine representation of $\V(M)$ inducing $\Omega_1(M)$.}
\end{rema}

\section{$H^i(\V(M),\cp^k_{diff}(\cc_M,\cs^1_2(M))), i=0,1.$}

We denote by $\cp^k_{diff}(\cc_M,\cs^1_2(M))$ the space of 
polynomials whose components relative to any base point are 
differential operators. It is clear that it is a $\V(M)$-submodule of 
$\cp^k(\cc_M,\cs^1_2(M))$.

\subsection{The main statement and its proof}

For the two first spaces of the corresponding cohomology, we have

\begin{theo}\label{calculcohomo}
\begin{eqnarray}
H^0(\V(M),\cp^k_{diff}(\cc_M,\cs^1_2(M)))&=&0, \ \forall k \geq 0, 
\label{un}\\
H^1(\V(M),\cp^k_{diff}(\cc_M,\cs^1_2(M)))&=&0, \ \forall k \geq 1, 
\label{deux}\\
H^1(\V(M),\cp^0_{diff}(\cc_M,\cs^1_2(M)))&\equiv& \rit^2.\label{trois}
\end{eqnarray}
\end{theo}

In this subsection, we describe the mains steps of the proof of that theorem, leaving the 
proofs of some auxiliary results for subsequent subsections.
The space  (\ref{trois}) has been computed in Proposition 
\ref{cohomologieS12}.
To prove (\ref{un}) and (\ref{deux}), we use the long exact sequence 
(\ref{longsequence}). It starts with
\begin{equation}\label{longsequence2}
\begin{array}{ccccccc}
0&\to& H^0(\cp^{k-1})&\to& H^0(\cp^k)&\to& H^0(S^k)\\
&\to& H^1(\cp^{k-1})&\to& H^1(\cp^k)&\to& H^1(S^k)\\
&\to&\cdots&&&&
\end{array}
\end{equation}
where we have simplified the notations by writing $H^i(\cp^l)$ instead of 
\[
H^i(\V(M),\cp^l_{diff}(\cc_M,\cs^1_2(M)))
\]
and $H^i(S^l)$ instead of
\[
H^i(\V(M),S^l_{diff}(\cs^1_2(M),\cs^1_2(M))).
\]
We are thus left to compute the latter, at least for $i=0,1$.
We will prove
\begin{prop}\label{auxiliaire}
\begin{eqnarray}
H^i(\V(M),S^k_{diff}(\cs^1_2(M),\cs^1_2(M)))&=&0, \ \forall i\in 
\{0,1\}, \forall k>1, \label{quatre}\\
H^1(\V(M),S^1_{diff}(\cs^1_2(M),\cs^1_2(M)))&\equiv&H^1(\V(M),C^\infty(M))\otimes\rit^2\label{huit}
\end{eqnarray}
\end{prop}
 Recall first \cite{14} that
\[
H^1(\V(M),C^\infty(M))=H^1_{DR}(M)\oplus\rit div_M
\]
where the subscript $DR$ denotes the de Rham cohomology of $M$ and $div_M$ is the {\it class of the divergence} of $M$, which is the class of the cocycle $X\mapsto D_\nabla(X)$ where $\nabla$ is the connection of Levi-Civita of some Riemannian metric of $M$.

Moreover, in (\ref{huit}), $\rit^2$ is spanned by $\tr.\1$ and $\pr$.
\begin{lema} \label{sequence}
For each $i\in \{0,1\}$ and each $k>1$, one has
\begin{eqnarray*}
H^i(\V(M),\cp^k(\cc_M,\cs^1_2(M)))=H^i(\V(M),\cp^1(\cc_M,\cs^1_2(M))).
\end{eqnarray*}
\end{lema}
\pre This follows immediately from the sequence (\ref{longsequence2}) and from (\ref{quatre}).$\cqfd$

\medskip
\begin{lema}\label{connectun}
The connecting homomorphism
\[
\chi:H^0(S^1)\to H^1(\cp^0)
\]
is an isomorphism.
\end{lema}
Using Proposition \ref{invariant} and Lemma \ref{connectun}, it follows immediately from the sequence (\ref{longsequence2}), with $k=1$, that $H^0(\cp^0)\equiv H^0(\cp^1)=0$, thus proving (\ref{un}), in view of Lemma \ref{sequence}.

\begin{prop}\label{connecdeux}
The connecting homomorphism
\[
\chi:H^1(S^1)\to H^2(\cp^0)
\]
is one-to-one.
\end{prop}

Equality (\ref{deux}) then follows immediately from Lemma \ref{sequence}, Lemma \ref{connectun} and from the above Proposition \ref{connecdeux}. Observe that (\ref{deux}) together with the exactness of the sequence (\ref{longsequence2}) implies Proposition \ref{connecdeux}. As we use that proposition to prove (\ref{deux}), we shall prove it directly. However, we are only able to do that when $M=\rit^m$. Therefore the proof of (\ref{deux}) will take two steps: first, we observe, as above, that it is true on $\rit^m$ using Proposition \ref{connecdeux}, second, we extend the result to arbitrary manifolds without using the above proposition, proving it in the same time on arbitrary manifolds.

We conclude this subsection by performing the second step.

We first observe that $1$-cocycle of $\V(M)$ valued in $\cp^k_{diff}(\cc_M,\cs^1_2(M))$ 
is local (this is easily seen, for instance like in \cite{12}).

We are now ready to prove (\ref{deux}) on an arbitrary manifold $M$ assuming that it is true on $\rit^m$.
Let $\gamma$ be a $1$-cocycle of $\V(M)$ valued in $\cp^1_{diff}(\cc_M,\cs^1_2(M))$, and  $\gamma_1$, $\gamma_0$ its homogeneous components of degree $1$ and $0$ relative to some $\nabla\in \cc_M$. Let $U$ an open subset of $M$ diffeomorphic to $\rit^m$. There exist $s_1^U\in\cd(\cs^1_2(U),\cs^1_2(U))$ and $s_0^U\in\cs^1_2(U)$ such that 
\[
\left\{
\begin{array}{rcl}
\gamma_{1\vert U}&=&\partial s_1^U\\
\gamma_{0\vert U}&=&\partial s_0^U-s_1^U\circ\cl_\nabla.
\end{array}
\right.
\]
We now prove that these are {\it unique} so that they are the restrictions over $U$ of some $s_1\in\cd(\cs^1_2(M),\cs^1_2(M))$ and $s_0\in\cs^1_2(M)$, showing that $\gamma$ is the coboundary of the polynomial whose components relative to $\nabla$ are $s_1,s_0$.

So, assume that $\partial s_1^U=0$ and that $\partial s_0^U-s_1^U\circ\cl_\nabla=0$.
In view of Proposition \ref{entrelacementS12}, it follows from the first equality that $s_1^U=a\tr.\1+b\pr$ for some real numbers $a,b$. The second equality then becomes 
\[
\partial s^U_0=a\cl_\nabla^{\tr}+b\cl_\nabla^{\pr}.
\] 
Since the cohomology classes of $\cl_\nabla^{\tr}$ and $\cl_\nabla^{\pr}$ are linearly independant, it follows that $a=b=0$. Hence $s^U_1=0$ and $\partial s^U_0=0$ so that $s^U_0=0$, due to Proposition \ref{invariant}. Hence the conclusion.

\subsection{Proof of Proposition \ref{auxiliaire}}

In order to compute the cohomology of $\V(M)$ acting on $S_{diff}^k(\cs^1_2(M),\cs^1_2(M))$, we filter that space by the order of differentiation. For each $r>0$, this gives a short exact sequence
\[
\begin{array}{l}
0\to S^{k,r-1}_{diff}(\cs^1_2(M),\cs^1_2(M))\to S^{k,r}_{diff}(\cs^1_2(M),\cs^1_2(M))\to\\[1ex]
P^{k,r}(\cs^1_2(M),\cs^1_2(M))\to 0
\end{array}
\]
where $S^{k,r}_{diff}(\cs^1_2(M),\cs^1_2(M))$ is the set of elements of $S_{diff}^k(\cs^1_2(M),\cs^1_2(M))$ of order at most $r$, $P^{k,r}(\cs^1_2(M),\cs^1_2(M))$ is the space of symbols of order $r$ and where the second map from the left is the inclusion.
This  sequence induces a long exact sequence of cohomology spaces, which begins this way:

\begin{equation}\label{longsequence3}
\begin{array}{ccccccc}0&\rightarrow& H^0(S^{k,r-1})&\rightarrow& 
H^0(S^{k,r})&\rightarrow& H^0(P^{k,r})\\
&\rightarrow& H^1(S^{k,r-1})&\rightarrow& H^1(S^{k,r})&\rightarrow& 
H^1(P^{k,r})\\&\rightarrow& \cdots&&&&
\end{array}
\end{equation}
(Again, we use shortcuts to avoid excessively lengthy notations.) 

\begin{prop}
For $(k,r)\neq (1,0)$,
\[
H^0(\V(M), P^{k,r}(\cs^1_2(M), \cs^1_2(M)))=0
\]
while
\[
H^0(\V(M), P^{1,0}(\cs^1_2(M), \cs^1_2(M)))=\rit^2.
\]
For $(k,r)\neq (0,0),(1,0)$,
\[
H^1(\V(M), P^{k,r}(\cs^1_2(M), \cs^1_2(M)))=0.
\]
Moreover
\[
H^1(\V(M), P^{0,0}(\cs^1_2(M), \cs^1_2(M)))=\rit \cl^\tr_M\oplus\rit\cl^\pr_M
\]
and
\[
H^1(\V(M), P^{1,0}(\cs^1_2(M), \cs^1_2(M)))=H^1(\V(M),C^\infty(M))\otimes\rit^2.
\]
\end{prop}

Again, $\rit^2$ is spanned by $\tr.\1$ and $\pr$.
The proof is quite long. It follows exactly the same line than the computation of $H^1(\V(M),S^1_{diff}(\cs^p(M),\cs^q(M)))$ exposed in \cite{13}. We will omit it.

It follows from that proposition and from (\ref{longsequence3}) that $H^0(S^{k,r})=H^0(S^{k,0})$ and $H^1(S^{k,r})=H^1(S^{k,0})$, for every $k,r\in\nit$. It remains to get rid of the order $r$ to conclude the proof of Proposition \ref{auxiliaire}.

We sketch the way to accomplish that. First, one shows that each $1$-cocycle of $\V(M)$ valued in 
$S^k_{diff}(\cs^1_2(M),\cs^1_2(M))$ is local as a map
\[
\V(M)\times\underbrace{\cs^1_2(M)\times...\times\cs^1_2(M)}_k\to\cs^1_2(M).
\]
Due to a theorem of Peetre, its restrictions on sufficiently small open sets of $M$ are differential operators, to which one can thus apply the above proposition. It remains then to glue the local informations so obtained to get the general form of the cocycle. Again, details will be omitted.

\subsection{Proofs of Lemma \ref{connectun} and Proposition \ref{connecdeux}}
The Lemma \ref{connectun} follows from the fact that the connecting homomorphism $\chi$ maps the basis $(\tr.\1,\pr)$ of $H^0(\cs^1)$ onto the basis $(\cl_M^{tr},\cl_M^{pr})$ of $H^1(\cp^0)$.

\medskip
In order to prove Proposition \ref{connecdeux}, we need two lemmas. Recall that the projective embedding $sl_{m+1}$ is the subalgebra of the algebra of vector fields with polynomial coefficients on $\rit^m$ generated by the fields $\partial_i, x^j\partial_i, x^ix^k\partial_k$ with
 $i,j=1,...,m$. It is isomorphic to $sl(m+1,\rit)$.
\begin{lema}\label{lemme1}
The space $H^2(sl_{m+1},\cs^1_2(\rit^m))$ is spanned by the class of the cocycle
\[
(X,Y)\mapsto \tr(\partial X)(L_Y\nabla^0)-\tr(\partial Y)(L_X\nabla^0),
\]
where $\nabla^0$ is the canonical flat connection on $\rit^m$ and $\partial X$ is the differential of the vector field $X$, given by $(\partial X)^i_j=\partial_j X^i$.
Moreover,
\[
H^1(sl_{m+1},\overline{\cs^1_2(\rit^m)})=0.
\]
\end{lema}
\pre This computation is made using a method  similar to the one described in \cite{13}. $\cqfd$

\begin{lema}\label{lemme2}
If the coboundary of $T\in Hom(\V(M),\cs^1_2(M))$ is local, then $T$ is local.
\end{lema}
\pre Let $U$ be an open subset of $M$ and $X$ a vector field on $M$ vanishing on $U$. Let $Y$ be any vector field with compact support included in $U$. We have then $T_{\lbrack X,Y\rbrack}=0$. Since $\partial T$ is local, this implies
\[
L_Y T_X|_U=0.
\]
As for the proof of Proposition \ref{invariant}, direct elementary computations in local coordinates show that $T_X|_U=0$. $\cqfd$ 

\medskip
Let us now suppose that $c\in H^1(\V(\rit^m), S_{diff}^1(\cs^1_2,\cs^1_2))$ is such that $\chi(c)=0.$
Proposition \ref{auxiliaire} provides us with two real $a$ and $b$ such that $c$ is the class of the cocycle $X\mapsto a\tr(\partial X)\tr.\1+b\tr(\partial X)\pr$. The cohomology class $\chi(c)$ is thus the class of the $2$-cocycle
\[
\begin{array}{rl} \kappa(X,Y)=& a\, \tr(\partial X)\tr(L_Y\nabla^0)\1+b\, \tr(\partial X)\pr(L_Y\nabla^0)\\
&-a\, \tr(\partial Y)\tr(L_X\nabla^0)\1-b\, \tr(\partial Y)\pr(L_X\nabla^0).
\end{array}
\]
As $\chi (c)$ vanishes, $\kappa$ is a coboundary.
Since $\pr(L_X\nabla^0)=0$ for $X\in sl_{m+1}$, we see that $\kappa (X,Y)=a(\tr(\partial X)(L_Y\nabla^0)-\tr(\partial Y)(L_X\nabla^0))$ if $X,Y\in sl_{m+1}$. Lemma \ref{lemme1} implies then $a=0$. We are now left to study the condition 
\begin{equation}\label{2cobord}
L_X T_Y-L_Y T_X-T_{\lbrack X,Y\rbrack}=b\,\tr(\partial X)\pr(L_Y\nabla^0))-b\,\tr(\partial Y)\pr(L_X\nabla^0)),
 \end{equation}
for some $T\in Hom(\V(\rit^m),\cs^1_2)$. We are going to show that $T$ is necessarily a $1$-cocycle. The proof of Proposition \ref{connecdeux} will then be achieved because the equality (\ref{2cobord}) becomes thus
\[
b\,\tr(\partial X)\pr(L_Y\nabla^0))-b\,\tr(\partial Y)\pr(L_X\nabla^0))=0,\;\forall X,Y\in\V(\rit^m),
\]
which obviously implies $b=0$.\\

In order to prove that $T$ is a $1$-cocycle, we are going to prove that the maps $\tr.\1\circ T$ and $\pr\circ T$ are $1$-cocycles.
Since the right hand side of equality (\ref{2cobord}) lies in $\overline{\cs^1_2}$, it is obvious that the map $\tr.\1\circ T$ is a $1$-cocycle. Substracting it from $T$, we can thus consider that $T$ is valued in $\overline{\cs^1_2}$. Let us now write down equality (\ref{2cobord}) for $X,Y\in sl_{m+1}$. As the right hand side of this equality then vanishes, $T_{\vert sl_{m+1}}$ is a $1$-cocycle. By Lemma \ref{lemme1}, $T$ is thus equal, up to a coboundary which can be extended to the whole $\V(\rit^m)$, to a map vanishing on $sl_{m+1}$.\\
We know by Lemma \ref{lemme2} and equation (\ref{2cobord}) that $T$ is local. Let us write (\ref{2cobord}) for some particular vector fields $X$, namely vector fields with polynomial coefficients of degre at most $1$ such that $\tr(\partial X)=0$. As $T_{\vert sl_{m+1}}=0$, we get for these $X$
\[
L_X T_Y-T_{\lbrack X,Y\rbrack}=0,\;\forall Y\in\V(\rit^m).
\]
Using the theory of invariants \cite{5}, we deduce that $T$ is necessarily a linear combination of $\cl_{\nabla^0}^{\tr}$ and $\cl_{\nabla^0}^{\pr}$, hence a constant multiple of the latter since $T$ is valued in $\overline{\cs^1_2}$. Therefore, $T$ is a $1$-cocycle.

\end{document}